\title{Boundary proximity of SLE}
\date{}
 \documentclass[12pt]{article}
   \usepackage{epsfig}
   \usepackage{color,latexsym,amsfonts,amssymb}
   \usepackage{amsmath}
\usepackage{amsthm}
\usepackage{graphicx}
\usepackage{hyperref}

\begin{document}
 \newcommand{\red}[1]{\textcolor{red}{#1}}
 \newcommand{\ignore}[1]{}{}
 \newcommand{\ul}[1]{{\bf #1}}
 \newcommand{\vv}{\vspace{8mm}}
 \newcommand{\vs}[1]{\vspace{#1cm}}
 \newcommand{\noi}{\noindent}
 \newcommand{\ds}{\displaystyle}
 \newcommand{\lra}{\longrightarrow}
 \newcommand{\la}{\lambda}
 \newcommand{\lrla}{\Longleftrightarrow}
 \newcommand{\beq}{\begin{eqnarray*}}
 \newcommand{\eeq}{\end{eqnarray*}}
 \newcommand{\beqn}{\begin{eqnarray}}
 \newcommand{\eeqn}{\end{eqnarray}}
 \newcommand{\ta}{\tau}
 \newcommand{\ra}{\rightarrow}
 \newcommand{\mb}{\mbox}
 \newcommand{\vp}{\varepsilon}
 \newcommand{\ep}{\epsilon}
 \newcommand{\al}{\alpha}
 \newcommand{\var}{\mbox{Var}}
 \newcommand{\bi}{\begin{itemize}}
 \newcommand{\ei}{\end{itemize}}
 \newcommand{\sta}{\stackrel}
 \newcommand{\La}{\Lambda}
 \newcommand{\Om}{\Omega}

 \newcommand{\ph}{\phi}
 \newcommand{\Ph}{\Phi}
 \newcommand{\ssb}{\scriptstyle \footnotesize 
                  \begin{array}{c}}
 \newcommand{\esb}{\end{array}}
 \newcommand{\PLDP}{PLDP }

 \newcommand{\nn}{\nonumber}
 \newcommand{\sq}{\sqrt}

 \newcommand{\lbl}{\label}
 \newcommand{\barB}{\bar{B}}
 \newcommand{\barb}{\bar{b}}
 \newcommand{\bars}{\bar{S}}
 \newcommand{\barv}{\bar{V}}
 \newcommand{\barx}{\bar{X}}
 \newcommand{\bbox}{\nobreak\quad\vrule width4pt depth2pt height4pt}
 \newcommand{\eq}[1]{$(\ref{#1})$}

\newtheorem{theorem}{Theorem}
\numberwithin{theorem}{section}
\newtheorem{corollary}[theorem]{Corollary}
\newtheorem{lemma}[theorem]{Lemma}
\newtheorem{remark}[theorem]{Remark}
\newtheorem{prop}[theorem]{Proposition}
\newtheorem{conjecture}[theorem]{Conjecture}
\newtheorem{problem}[theorem]{Problem}
\newtheorem{exercise}[theorem]{Exercise}

 \renewcommand{\theequation}{\arabic{section}.\arabic{equation}}

 \newcommand{\f}{\frac}
 \newcommand{\wt}{\widetilde}
 \renewcommand{\(}{\left(}
 \renewcommand{\)}{\right)}

 \newcommand{\sn}{\sum_{i=1}^{n_2}}

 \renewcommand{\P}{{\bf P}}

 \newcommand{\E}{{\bf E}}
 \newcommand{\B}{{\bf B}}

 \newcommand{\de}{\delta}
 \newcommand{\De}{\Delta}
 \newcommand{\ga}{\gamma}
 \newcommand{\Ga}{\Gamma}
 \newcommand{\si}{\sigma}
 \newcommand{\Si}{\Sigma}
 \newcommand{\Th}{\Theta}
 \newcommand{\pa}{\parallel}
 \newcommand{\ka}{\kappa}
 \newcommand{\no}{\nonumber}
 \newcommand{\ol}{\overline}

\newif\ifdraft
\drafttrue \long\def\note#1{\ifdraft{\marginpar{{$\Longleftarrow$}} \bf [#1] }\fi}
 \def\eref#1{(\ref{#1})}
\def\dhrule{\bigskip\hrule\smallskip\hrule\bigskip}
\def\QED{\qed\medskip}
\newcommand{\Prob} {{\bf P}}
\newcommand{\R}{\mathbb{R}}
\newcommand{\C}{\mathbb{C}}
\newcommand{\Z}{\mathbb{Z}}
\newcommand{\N}{\mathbb{N}}
\newcommand{\HH}{\mathbb{H}}
\def\H{\mathbb{H}}
\def\U{\mathbb{U}}
\def\diam{\mathrm{diam}}
\def\area{\mathrm{area}}
\def\dist{\mathrm{dist}}
\def\ceil#1{\lceil{#1}\rceil}
\def\floor#1{\lfloor{#1}\rfloor}
\def\length{\mathrm{length}}
\def\Re{\operatorname{Re}}
\def\Im{\operatorname{Im}}
\def\SLEkk#1/{$\mathrm{SLE}_{#1}$}
\def\SLEk/{\SLEkk{\kappa}/}
\def\SLE{\mathrm{SLE}}
\def\Ito/{It\^o}
\def \eps {\epsilon}
\def \P {{\bf P}}
\def\md{\mid}
\def\Bb#1#2{{\def\md{\bigm| }#1\bigl[#2\bigr]}}
\def\BB#1#2{{\def\md{\Bigm| }#1\Bigl[#2\Bigr]}}
\def\Bs#1#2{{\def\md{\mid}#1[#2]}}
\def\Pb{\Bb\P}
\def\Eb{\Bb\E}
\def\PB{\BB\P}
\def\EB{\BB\E}
\def\Ps{\Bs\P}
\def\Es{\Bs\E}
\def \prob {{\bf P}}
\def \expect {{\bf E}}
\def \p {{\partial}}
\def \E {{\bf E}}
\def\defeq{:=}
\def\closure{\overline}
\def\ev#1{{\mathcal{#1}}}
\def \proof {{ \medbreak \noindent {\bf Proof.} }}
\def\proofof#1{{ \medbreak \noindent {\bf Proof of #1.} }}
\def\proofcont#1{{ \medbreak \noindent {\bf Proof of #1, continued.} }}

\def\zzz{\bigskip\hrule\smallskip\hrule\bigskip}
\def\xo{\xi}
\def\field{\mathcal{F}}
\def\ix{\hat x}
\def\iu{\hat u}
\def\iv{\hat v}
\def\iy{\hat y}
\def\iz{\hat z}
\def\tu{\tilde u}
\def\hxo{\hat\xo}
\def\hf{\hat f}
\def\hh{\vartheta}
\def\GG{G}
\def\hG{\hat G}
\def\tG{\widetilde G}
\def\FF{\tilde F}
\def\hF{\hat F}
\def\bF{\bar F}
\def\Li{\mathrm{Li}}
\def\dd{d}
\def\sb{{s_\kappa}}
\def\Hk{\Lambda_\kappa}
\def\QQ{Q_a}
\def\yy{h}
\def\aa{a}
\def\bb{b}

\def\sign{\mathrm{sign}}
\def\BCM{MR1815718}
\def\BMakFine{MR2000g:30018}
\def\BJonMak{MR96k:30027}
\def\BDurStoch{MR97k:60148}
\def\BSmirnovPerc{MR1851632}
\def\BDup{MR2001c:82040}
\def\BDupSal{MR88d:82073}
\def\WernerSurvey{MR1905353}
\def\BhigherFunc{MR15:419i}
\def\BCardySurvey{math-ph/0103018}
\def\BLSWi{MR2002m:60159a}
\def\BLSWii{MR2002m:60159b}
\def\BLSWiii{MR1899232}
\def\BLSWlesl{math.PR/0112234}
\def\BLSWrest{MR1992830}
\def\BSchSLE{MR1776084}
\def\BPomm{MR95b:30008}
\def\BAhlfors{MR50:10211}
\def\BKR{MR98k:28004}
\def\BMR{MRlow}
\def\BRHAdd{MR98e:82018}
\def\BWeilWil{WeilWil}
\def\BBeffaraHaus{math.PR/0211322}
\def\DubedatDuality{math.PR/0303128}
\def\RevuzYor{MR2000h:60050}
\def\SSharmonicExplorer{math.PR/0310210}

\author{%
Oded SCHRAMM\thanks{Microsoft Research}
\and
Wang ZHOU\thanks{National University of Singapore; supported in part by grant
R-155-050-055-133/101 at the National University of Singapore}
}
\maketitle

\begin{abstract}
This paper examines how close the chordal $\SLE_\kappa$ curve gets to the real line
asymptotically far away from its starting point.
In particular, when $\kappa\in(0,4)$, it is shown that if
$\beta>\beta_\kappa:=1/(8/\kappa-2)$, then the intersection of the $\SLE_\kappa$ curve
with the graph of
the function $y=x/(\log x)^{\beta}$, $x>e$, is a.s.\ bounded,
while it is a.s.\ unbounded if $\beta=\beta_\kappa$.
The critical $\SLE_4$ curve a.s.\ intersects the graph of
$y=x^{-(\log\log x)^\alpha}$, $x>e^e$, in an unbounded set
if $\alpha\le 1$, but not if $\alpha>1$.
Under a very mild regularity assumption on the function $y(x)$,
we give a necessary and sufficient integrability condition for the intersection
of the $\SLE_\kappa$ path with the graph of $y$ to be unbounded.
When the intersection is bounded a.s.\,, we provide an estimate for
the probability that the $\SLE_\kappa$ path hits the graph of $y$.
 We also prove that  the Hausdorff dimension of the intersection set of the $\SLE_{\kappa}$ curve
and real axis is $2-8/\kappa$ when $4<\kappa<8$.
\end{abstract}



\bigskip

\section{Introduction}
The stochastic Loewner evolution paths ($\SLE$) are random curves in the plane that
are obtained by running Loewner's differential equation with a scaled Brownian motion
as the driving parameter. They have been shown to describe several critical
statistical physics systems, and have been useful in the analysis of these systems.
This has been proved for critical site-percolation on the triangular lattice
\cite{smirnov01,CN}, loop erased random walks and uniform spanning tree Peano paths
\cite{LSW04}, the level lines of the discrete Gaussian free field~\cite{SS08},
the interfaces of the random cluster model associated with the Ising model~\cite{smirnov07},
as well as a few other systems.
For further background, the reader is advised to consult the
surveys~\cite{werner,GK,cardy,lawler}.

In order to understand the corresponding disordered systems well,
it is then natural to investigate the properties of $\SLE$.
In \cite{RS}, the basic topological and geometric properties of $\SLE$
were investigated.
In \cite{beffara04}
the Hausdorff dimension of the $\SLE_6$ curve and its outer boundary
were determined.
Several years later it was proved~\cite{beffara} that the
Hausdorff dimension of the $\SLE$ curves is $\min(1+\kappa/8,2)$.

There are several different versions of $\SLE$.
If $B_t$ is a one-dimensional Brownian motion starting at $0$,
the chordal $\SLE_{\kappa}$ in the upper half plane $\H$ from $0$ to $\infty$ with parameter
$\kappa$ is the solution of the differential equation
\beqn
\partial_t g_t(z)= \f{2}{g_t(z)-W_t}, \ \ g_0(z)=z,  \lbl{chordal}
\eeqn
where $z\in \H$ and $W_t=\sqrt{\kappa}B_t$.
It can be shown~\cite{RS, LSW04} that a.s.\  $g_t^{-1}$ extends continuously
to $\closure\H$ for every $t\ge 0$ and $\gamma(t):=g_t^{-1}(W_t)$ is a continuous path.
This is the $\SLE$ path, and the domain of definition of $g_t$ is the unbounded
connected component $H_t$ of $\H\setminus\gamma[0,t]$.
We shall denote by $K_t$ the closure of the complement of $H_t$ in $\H$.

It is known~\cite{RS} that when $\kappa\ge 8$ a.s.\ $\gamma\cap \R=\R$ and
when $\kappa\in[0,4]$ a.s.\ $\gamma\cap \R=\{0\}$.
In this paper, we will study the boundary behavior of $\SLE$ curves.
More precisely, given the graph of a function $\yy:[r,\infty)\to(0,\infty)$
we will discuss whether the intersection set of the $\SLE_{\kappa}$ curve $\gamma$ and
the graph of $\yy(x)$ is bounded or not.
Clearly, this intersection is a.s.\ unbounded
 when $\kappa>4$, since $\gamma$ swallows every point of $\H$ a.s.\ when $4<\kappa<8$ and
$\gamma=\closure{\H}$ when $\kappa\geq 8$. The only non-trivial case is
$\kappa \in(0,4]$.

For a function $\yy:[r,\infty)\to(-\infty,\infty)$, let $\Gamma^\yy$ denote its graph; that is,
$$
\Gamma^\yy:=\{x+i\,\yy(x):x\ge r\}\,.
$$
Set
$$
\sb :=8/\kappa-1\,,
$$
and
$$
\Hk^h(x):=\begin{cases} h(x)^{\sb -1}&\kappa<4\,,\\
             1/\log\big(\frac x{h(x)}\vee 2\big)&\kappa=4\,.
\end{cases}
$$
Our main theorem is the following.

\begin{theorem}\label{t.main}
Let $\kappa\in(0,4]$, and let $\gamma$ be the chordal $\SLE_\kappa$ path.
Fix $r>1$, and suppose that $\yy:[r,\infty)\to(0,\infty)$ is continuous
and satisfies
\begin{equation}
\label{e.regularity}
\sup\bigl\{\Hk^h(x)/\Hk^h(y):r\le x\le y\le 2\,x\}<\infty\,.
\end{equation}
If
\begin{equation}
\label{e.test} \int_r^\infty \frac {\Hk^{h}(x)}{x^{\sb
}}\,dx<\infty \,,
\end{equation}
then $\gamma\cap\Gamma^\yy$ is bounded a.s.
Conversely, if the integral in~\eqref{e.test} is infinite,
then $\gamma\cap\Gamma^\yy$ is unbounded a.s.
\end{theorem}
To illustrate the theorem, we note that if $\kappa<4$ and $\yy(x)=x\,(\log x)^{-\beta}$,
then $\gamma\cap\Gamma^\yy$ is bounded a.s.\ if $\beta>(8/\kappa-2)^{-1}$
and unbounded a.s.\ if $\beta=(8/\kappa-2)^{-1}$.

The case $\kappa=4$ is critical for $\SLE$ to hit the boundary, and it is therefore
not entirely surprising that its behavior is different. In that case, if
$\yy(x)=x^{-(\log\log x)^\alpha}$, then $\gamma\cap\Gamma^\yy$ is a.s.\
unbounded when $\alpha=1$, but bounded a.s.\ if $\alpha>1$.

Now suppose instead that $\yy$ is continuous in $[0,1]$ and $\yy(0)=0$.
One can ask if $0$ is in the closure of the intersection of $\bigl\{x+i\,h(x):x\in (0,1]\bigr\}$
and $\gamma$. Using reversibility of SLE~\cite{zhan07}, this translates
to the type of question addressed by Theorem~\ref{t.main}.
Alternatively, the proof of Theorem~\ref{t.main} can be easily adapted to also handle
this question.

Another natural question related to Theorem~\ref{t.main} is to estimate
the probability that $\gamma$ hits $\Gamma^\yy$.
Actually, it is not too hard to modify the proof of Theorem~\ref{t.main} to show
that when $\kappa\le 4$
\begin{equation}
P(\gamma\cap \Gamma^\yy\ne\emptyset)
\eqsim 1\wedge\int_r^\infty
\frac{\Hk^{h}(x)}{x^{\sb }}\,dx\,,
\end{equation}
where $\eqsim$ denotes equivalence
up to a multiplicative constant that
depends only on $\kappa$ and the left hand side in~\eqref{e.regularity}.
Likewise, the proof of Theorem~\ref{t.main} easily gives the following
estimate for the probability that $\gamma$ hits the
set $A_{\epsilon}=\{x+i\,y:1\leq x\leq 2,\,0\le y\le\eps\}$:
\begin{equation}
P(\gamma \cap A_{\epsilon}\not=\emptyset)\asymp
\begin{cases} \epsilon^{\sb-1} &\kappa<4\,,\\
        |\log \epsilon|^{-1} &\kappa=4\,,
\end{cases}
\end{equation}
where $\eps\in(0,1/2)$ and
the constants implied by $\asymp$ depend only on $\kappa$.
Somewhat related results in the setting of discrete models appear
in~\cite[Theorem 10.7]{sch00} and in~\cite{vdB-Jarai}.

\medskip

We also make use of the machinery developed for the proof of
Theorem~\ref{t.main}
to obtain the Hausdorff dimension of $\gamma\cap \R$ when $\kappa\in(4,8)$.

\begin{theorem} \label{theore3}
If $4<\kappa<8$, then with probability one,
$$
\dim_H \bigl(\gamma \cap \R\bigr)=2-8/\kappa\,.
$$
\end{theorem}

A proof of this result based on Beffara's argument should be possible,
but our proof is different and simpler. In fact, one may hope that
the argument we present would generalize to give a simpler proof of
Beffara's theorem, but so far we were not able to achieve this.
An alternative and independent proof of Theorem~\ref{theore3}
can be found in~\cite{AS}.
\medskip

The paper is organized as follows. In Section~\ref{mm}, we consider
for each $x>0$ a local martingale $M^x_t$ and relate its behavior
to the geometry of the path near $x$.
We also derive an estimate for the probability that both
$M^x_t$ and $M^y_t$ become large, as a function of the
positions of the points $x,y$.
In Section~\ref{G}, we prove Theorem~\ref{t.main}
 using the first and second moment methods.
The Hausdorff dimension proof is given in Section~\ref{HD}.

\section{The local martingale and its properties}\label{mm}
\subsection{Basic properties}\label{ub}

We assume throughout this paper that $\kappa\in(0,8)$.
Let $x>0$ and set
$$
t_x:= \sup\{t\ge 0: x\notin K_t\}\,.
$$
Then we have from~\cite{RS} that $t_x=\infty$ a.s.\ if $\kappa\le 4$ and $t_x<\infty$ a.s.\
if $\kappa>4$.
Define for $t\in(0,t_x)$,
$$
M^x_t:=\Bigl(\frac{g'_t(x)}{g_t(x)-W_t}\Bigr)^{\sb }\,.
$$
Also, for $\eps>0$ set
$$
\tau_x=\tau_x^\eps:=\inf\bigl\{t\in(0,t_x):M^x_t\ge \eps^{-\sb }\bigr\}
$$
and
\begin{equation} \label{cepsilon}
C_\eps:=\{x>0:\tau_x<t_x\}\,.
\end{equation}
As usual, we use the convention that $\inf\emptyset=\infty$.

Write $\mathcal{F}_t:=\sigma(B_s, 0\leq s\leq t)$. Then
$\{M^x_t, \mathcal{F}_t, t\in (0,t_x)\}$ is a local martingale by Theorem 6 and Remark 7
in \cite{SW} (this is, of course, easily verified using It\^o's formula).
The reason for our interest in $M_t^x$ is the following lemma.

\begin{lemma}\label{l.dist}
If $x\in C_\eps$, then the distance from $\gamma$ to $x$ is at most $4\,\eps$.
\end{lemma}
\proof
Suppose that $x>0$, $t>0$, and $x\notin K_t$.  Set $\bar K_t:=\{\bar z:z\in K_t\}$, and let
$G$ denote the extension of $g_t$ to $\C\setminus (K_t\cup\bar K_t)$,
which is obtained by Schwarz reflection.
Let $\dd_t=\dd^x_t$ denote the distance from $x$ to $\gamma[0,t]$.
Then $W_t$ is not in $G\bigl(B(x,\dd_t)\bigr)$, and therefore the Koebe $1/4$ theorem
gives
$$
G'(x)\,\dd_t/4\le G(x)-W_{t}\,.
$$
This translates to
\begin{equation}
 \label{e.dist}
M_t^x \le (4/\dd_t)^{\sb },
\end{equation}
and the lemma immediately follows.
\QED

Next, we prove that in some situations the inequality~\eqref{e.dist} may be reversed.

\begin{lemma}\label{l.rev}
Let $x>0$, $t>0$, $x_0:=\Re \gamma(t)$ and $y_0:=\Im\gamma(t)$.
Suppose that $x\notin K_t$, $x-x_0\ge y_0$,
and $\gamma[0,t)$ does not intersect the line segment $[x_0,\gamma(t)]$.
Then
$$
M_t^x \ge (c\,\dd_t)^{-\sb },
$$
where $0<c<\infty$ is a universal constant.
\end{lemma}

\proof
Let $G$ be as in the proof of Lemma~\ref{l.dist}.
Set $r:=G(x)-G(x_0)$.
Then the inverse of $G$ is defined in the ball $B\bigl(G(x),r\bigr)$.
Therefore, the Koebe $1/4$ theorem gives
$$
\frac r 4\,G'(x)^{-1}\le \dd_t\,.
$$
It therefore suffices to prove a positive lower bound on
\begin{equation}
\label{e.bb}
\frac{r}{G(x)-W_t}=
\frac{G(x)-G(x_0)}{G(x)-W_t}\,.
\end{equation}
Every path in $H_t$ going from $[x,\infty)$ to the union of $[0,x_0]$ and
the right hand side of $\gamma[0,t]$
must intersect the line segment $[x_0,\gamma(t)]$.
Since we may consider the Euclidean metric on
the square of sidelength $2\,y_0$ centered at $\gamma(t)$,
normalized to have area $1$,
the extremal length of this
collection of paths is bounded away from zero.
By conformal invariance of extremal distance,
it follows that the extremal distance from $[G(x),\infty)$ to $[W_t,G(x_0)]$ in $\H$
is likewise bounded away from zero.
This implies the required lower bound on~\eqref{e.bb}, and completes the proof.
\QED

For a given point $x>0$, we are interested in the probability
that $x\in C_{\epsilon}$.

\begin{prop} \label{m1}
Let $0<\kappa<8$, $x>0$ and $\eps>0$. Then
\begin{equation}
\label{e.m1}
    P(x \in C_{\epsilon} )
=(\epsilon/x)^{\sb }\wedge 1\,.
\end{equation}
\end{prop}
The proof is dependent on the properties of the local martingale
$M^x_t$ as $t\nearrow t_x$.
Write $T_x=t_x\wedge \tau_x$.
If $0<\kappa\leq 4$, then $t_x=\infty$ a.s.\ and $T_x=\tau_x$.
We use $I(\ev A)$ for the indicator function of an event $\ev A$.

\proof
Since $M_{t\wedge T_x}^x$ is a bounded local martingale on $t\in(0,T_x)$,
it is a martingale and the limit $M_{T_x}^x:=\lim_{t\nearrow T_x} M^x_t$ exists.
On the event $0<\tau_x<\infty$, we have $M_{T_x}=\eps^{-\sb }$.
Hence, the optional sampling theorem gives
$$
M_0^x =
E(M_{T_x}^x)= P(x\in C_\eps)\, \eps^{-\sb } + E\bigl( M^x_{T_x}\,I(\tau_x=\infty)\bigr)\,.
$$
Therefore, the proof is complete once we prove that
\begin{equation}\label{e.swallow}
P(M_{T_x}^x\ne 0,\,\tau_x=\infty)=0\,.
\end{equation}

Consider first the case $\kappa\in(4,8)$.
In this case a.s.\ $t_x<\infty$ and $x_1:=\gamma(t_x)\in(x,\infty)$.
Suppose that this is indeed the case.
Let $r>0$ be much smaller than the distance from $x$ to $x_1$, and let
$s$ be the first time $t$ at which $|\gamma(t)-x_1|=r$.
Let $G$ denote the Schwarz reflection of $g_s$ with respect to the real line,
let $a=a_s:=\sup (K_s\cap\R)$ and $a'=a'_s:=\inf\{G(x'):x'>a\}$.
Then $a'$ is not in $G\bigl(B(x,\dd_s)\bigr)$. Therefore, the Koebe $1/4$ theorem
implies
$$
G'(x)\,\dd_s/4 \le  G(x)-a'\,.
$$
That is,
$$
\frac{g_s'(x)}{g_s(x)-a'}\le 4/\dd_s\,.
$$
Since $\dd_{t_x}>0$ a.s.,
it therefore suffices to prove that
\begin{equation}
\label{e.rat}
\lim_{r\searrow 0}\; \frac{g_s(x)-W_s}{g_s(x)-a'}=\infty\,.
\end{equation}

Consider the extremal distance in $H_s$ from $(a,x)$ to the union of
$(-\infty,0)$ and the left hand side of $\gamma[0,s]$.
This extremal distance is clearly at least as large as the extremal distance from the circle
of radius $|x-x_1|$ about $x_1$ to the circle of radius $r$ about $x_1$, which is
at least a constant times $\log \bigl(|x-x_1|/r\bigr)$.
By conformal invariance of extremal distance, it follows that the
extremal distance in $\H$ from $[a',g_s(x)]$ to $(-\infty,W_s]$
goes to infinity as $r\searrow0$, which proves~\eqref{e.rat} and
completes the proof in the case $\kappa\in(4,8)$.

The argument in the case $\kappa\in(0,4]$ is similar.
We choose $R>0$ large, and let $s$ be the first time at which $|\gamma(s)|=R$.
The extremal distance in $H_s$ from $(0,x]$ to the union of the left hand side of
$\gamma[0,s]$ with $(-\infty,0)$ is then at least a constant times $\log(R/x)$,
which implies~\eqref{e.rat} in the same way.
\QED

Observe that the proposition implies that given $x>0$ there is a.s.\ some $\eps>0$
such that $x\notin C_\eps$.  Therefore,~\eqref{e.swallow} gives
\begin{equation}
\label{e.Mlim}
M_{t_x}^x:= \lim_{t\nearrow t_x} M_t^x=0\qquad \text{a.s.}
\end{equation}

\subsection{Correlation estimate} \label{lb}
Let $0<x<y$, $\eps_x,\eps_y>0$, $\tau_x:=\tau_x^{\eps_x}$,
$\tau_y:=\tau_y^{\eps_y}$,
$T_x:=t_x\wedge\tau_x$ and $T_y:=t_y\wedge\tau_y$. Define
$$Z_t:=\frac{g_t(x)-W_t}{g_t(y)-W_t}, \ \ T:=T_x\wedge T_y\,.$$
A simple but tedious calculation via It\^o's formula implies that $u(Z_t)\,M^x_t\,M^y_t$ is
a local martingale while $t\in(0,T)$, where
$$
u(z):=\,(1-z)^{-\sb } \,{}_2F_1(1-8/\kappa,4/\kappa,8/\kappa;1-z)\,.
$$
 Euler's integral representation of hypergeometric functions shows that
\begin{multline*}
{}_2F_1(1-8/\kappa,4/\kappa,8/\kappa;z)
=\frac{\Gamma(8/\kappa)}{\Gamma(4/\kappa)^2}\int_0^1t^{4/\kappa-1}
(1-t)^{4/\kappa-1}(1-zt)^{8/\kappa-1}\,dt,
\end{multline*}
where $\Gamma$ is the gamma function.
This implies that $u(z)>0$ when
$z\in (0,1)$.
Since $8/\kappa-4/\kappa-(1-8/\kappa)>0$, we have
$$
{}_2F_1(1-8/\kappa,4/\kappa,8/\kappa;1)=\frac{\Gamma(8/\kappa)\Gamma(12/\kappa-1)}
{\Gamma(16/\kappa-1)\Gamma(4/\kappa)}\,.
$$
Hence $q_1:=\inf_{z\in(0,1)}u(z)$
and $q_2:=\sup_{z\in(0,1)} (1-z)^{\sb }\,u(z)$ are both finite and positive.
It follows from~\eqref{e.swallow} that
\begin{equation*}
P(x\in C_{\epsilon_x}, y\in C_{\epsilon_y})=P(M_{T_x}^x=\eps_x^{-\sb },\,M_{T_y}^y=\eps_y^{-\sb })
=(\eps_x\,\eps_y)^{\sb }\,E(M_{T_x}^x\,M_{T_y}^y).
\end{equation*}
Recall that $T=T_x\wedge T_y$.
If $T=T_x<\infty$, then we have that $M^x_{t\wedge T_x}$ is constant
in the range $t\in[T_x,T_y)$, while $M^y_{t\wedge T_y}$ is a martingale.
The symmetric statement also holds when we exchange $x$ and $y$.
It should be clear that this implies
\begin{equation}
\label{e.marteq}
E(M_{T_x}^x\,M_{T_y}^y)=E(M_{T}^x\,M_{T}^y)\,,
\end{equation}
but for the sake of completeness, we prove this.
First, since $I(T=T_x)\,M_{T_x}^x$ is $\mathcal F_T$-measurable, we have
$$
\begin{aligned}
E\bigl(I(T=T_x)\,M_{T_x}^x\,M_{T_y}^y\bigm| \mathcal F_{T}\bigr)
&
=
I(T=T_x)\,M_{T_x}^x\,E\bigl(M_{T_y}^y\bigm| \mathcal F_{T}\bigr)
\\&
= I(T=T_x)\,M_T^x \,M_T^y\,.
\end{aligned}
$$
Second, on the complement of the event $T=T_x$, we have $T=T_y$. Hence,
$I(T\ne T_x)\,M_{T_y}^y$ is also $\mathcal F_T$-measurable, and we get in the same way
$$
E\bigl(I(T\ne T_x)\,M_{T_x}^x\,M_{T_y}^y\bigm| \mathcal F_{T}\bigr)
= I(T\ne T_x)\,M_T^x \,M_T^y\,.
$$
Summing the above and taking expectations, we obtain~\eqref{e.marteq}.

Since $u(Z_{t\wedge T})\,M_{t\wedge T}^x\,M^y_{t\wedge T}$ is a non-negative local martingale,
it is also a super\-martingale. This justifies the second inequality in the
following estimate:
\begin{equation*}
\begin{aligned}
P(x\in C_{\epsilon_x}, y\in C_{\epsilon_y})&=
(\eps_x\,\eps_y)^{\sb }\, E(M^x_T\,M^y_T)
\\&
\le
(\eps_x\,\eps_y)^{\sb }
\,E\bigl(u(Z_T)\, M^x_T\,M^y_T\bigr)/q_1
 \\ &
\le
(\eps_x\,\eps_y)^{\sb }
\,u(Z_0)\, M^x_0\,M^y_0/q_1
 \\ &
=
(\eps_x\,\eps_y)^{\sb }
\,\frac{u(x/y)}{q_1\,x^{\sb }\,y^{\sb }}
\le \frac {q_2}{q_1}\,
(\eps_x\,\eps_y)^{\sb }
\,x^{-\sb }\,(y-x)^{-\sb }\,.
\end{aligned}
\end{equation*}
Hence, we obtain the following proposition.

\begin{prop} \label{m2}
Let $0<x<y$, $0<\kappa<8$, $\eps_x,\eps_y>0$. Then
\begin{equation}
\label{e.m2}
P(x\in C_{\epsilon_x}, y\in C_{\epsilon_y})
\leq c_{\kappa}\,(\eps_x\,\eps_y)^{\sb }\,x^{-\sb }\,(y-x)^{-\sb },
\end{equation}
where $c_{\kappa}$ is a constant depending only on $\kappa$. \QED
\end{prop}

\section{Proximity estimates} \label{G}
\subsection{Bounded intersection} \label{s.lb}

In this subsection, we assume~\eqref{e.test}, as well
as the other assumptions in Theorem~\ref{t.main},
and prove that $\gamma\cap\Gamma^\yy$ is bounded a.s.

\medskip

Let $\rho:[r,\infty)\to(0,\infty)$ be a function such that
\begin{equation}
\label{e.rhoup}
\lim_{x\to\infty}\rho(x)/\Hk^h(x)=\infty\,,
\end{equation}
but $\rho$ satisfies~\eqref{e.regularity} and~\eqref{e.test}
in place of $\Hk^h(x)$, namely,
\begin{equation}
\label{e.rhodef}
\int_r^\infty \frac {\rho(x)}{x^{\sb }}\,dx  <\infty\,,
 \end{equation}
and
 \begin{equation}
\label{e.rhobd}
\sup\bigl\{\rho(x)/\rho(y):r\le x\le y\le 2\,x\bigr\}  <\infty\,.
 \end{equation}
In the following, we let $\eqsim$ mean equivalence up to
positive multiplicative
constants that may depend on $\yy, \rho$ and $\kappa$. Likewise,
$a \lesssim b$ means that there is some $a'\le b$ such that $a'\eqsim a$.

Define
$$
Z_t:=\int_r^\infty \rho(x)\,M^x_t\,dx\,.
$$
Since $M^x_0=x^{-\sb }$, it follows from~\eqref{e.rhodef} that $Z_0<\infty$.
As $M^x_t$ is a super\-martingale for each $x>0$, it follows that
$Z_t$ is a super\-martingale.

By~\eqref{e.regularity}, for every $x\ge r$ such that $\yy(x)\ge x/2$,
the contribution to the integral in~\eqref{e.test} from the interval $[x,2\,x]$
is bounded from zero.
Since the integral in~\eqref{e.test} is finite, we conclude that there is
a finite $R_0>r$ such that $\yy(x)<x/2$ for $x\ge R_0$.
Fix an $R>R_0$, and let $A$ be the set $A:=\{x+i\,y:x\ge R,\,y\le \yy(x)\}$.
Let $T_A:=\inf\{t\ge 0:\gamma_t\in A\}$, and on the event $T_A<\infty$
set $x_0:=\Re \gamma(T_A)$, $y_0:=\Im \gamma(T_A)$. From our choice of $R$,
we have $y_0\le x_0/2$.
By Lemma~\ref{l.rev}, $M_{T_A}^x \gtrsim (x-x_0)^{-\sb }$
holds for every $x>x_0+y_0$.
Therefore, on the event $T_A<\infty$,
\begin{multline*}
Z_{T_A}
\gtrsim \int_{x_0+y_0}^\infty \rho(x)\,(x-x_0)^{-\sb }\,dx
\overset{\eqref{e.rhobd}}
\gtrsim \rho(x_0)\int_{x_0+y_0}^{2 x_0} (x-x_0)^{-\sb }\,dx
\\
\gtrsim
\frac{\rho(x_0)}{\Hk^h(x_0)}
\,.
\end{multline*}
Since $Z_t$ is a supermartingale,  the optional sampling theorem
gives
$$
Z_0\ge E\bigl(Z_{T_A}\,I(T_A<\infty)\bigr)
\gtrsim P(T_A<\infty)\,\inf_{x\ge R}
\frac{\rho(x)}{\Hk^h(x)}
\,.
$$
By~\eqref{e.rhoup}, we conclude that
$\lim_{R\to\infty} P(T_A<\infty)=0$.
Thus, $\gamma\cap\Gamma^\yy$ is bounded a.s., as required. \QED

\subsection{Unbounded intersection} \label{UI}

In this subsection, we assume that the integral in~\eqref{e.test}
is infinite, and prove that $\gamma\cap\Gamma^\yy$ is unbounded a.s.,
thus completing the proof of Theorem~\ref{t.main}.
In the following, $\eqsim$ denotes equivalence up to multiplicative
constants that may depend on $\kappa$ and $\yy$,
and similarly for $\lesssim$.

\medskip

Suppose that we prove that the intersection of
$\gamma$ with
$\Theta^h_+:=\bigl\{x+i\,y:y\le h(x),\,x\ge r\bigr\}$ is  a.s.\
unbounded. Symmetry then implies that the
intersection of $\gamma$ with
$\Theta^h_-:=\bigl\{-x+i\,y:y\le h(x),\,x\ge r\bigr\}$ is a.s.\ unbounded
as well.
For every $R\ge \sqrt{r^2+h(r)^2}$,
the set $\Gamma^h\setminus B(0,R)$ separates
$\Theta^h_+\setminus B(0,R)$ from
$\Theta^h_-\setminus B(0,R)$ in $\H\setminus B(0,R)$.
Since $\gamma$ is a.s.\ transient,
it follows that $\gamma\cap \Gamma^h$ must be a.s.\ unbounded, as required.

Now observe that $\tilde h(x):=h(x)\wedge (x/2)$ satisfies the same
assumptions as we have for $h$. The above then implies that it suffices
to prove the claim for $\tilde h$.
Thus, we assume with no loss of generality that $h(x)\le x/2$
holds for every $x\ge r$.

Define $\rho(x):=\Hk^{h}(x)$.
Let $\aa>r$ and let $\bb>\aa$ satisfy
\begin{equation}
\label{e.Rp}
\int_{\aa}^{\bb} \frac{\rho(x)}{x^{\sb}}\,dx=1\,.
\end{equation}
Define $X:=\{x\ge r:x\in C_{\yy(x)}\}$.
We will show that $\sup X=\infty$ a.s.
Set
$$
\QQ:=\int_{\aa}^{\bb} \frac{\rho(x)}{\yy(x)^\sb}\,I(x\in X)\,dx\,.
$$
Then by~\eqref{e.Rp} and Proposition~\ref{m1}, we have
$$
E(\QQ)=1\,.
$$

We will now prove that $E(\QQ^2)$ is bounded by some constant independent
of $\aa$. First, observe that
$$
E(\QQ^2)=
\int_{\aa}^{\bb}
\int_{\aa}^{\bb}
\frac{\rho(x)\rho(y)}{\yy(x)^\sb \yy(y)^\sb}\,P(x,y\in X)\,dx\,dy
\,.
$$
Let $F(x,y)$ denote the integrand. Set $S:=[\aa,\bb]^2$.  Let
$S_1$ be the set of pairs $(x,y)\in S$ such that $y\in[x-\yy(x),x]$,
let $S_2$ be the set of pairs $(x,y)\in S$ such that $y\in[x/2,x-\yy(x)]$,
and let $S_3$ be the set of pairs $(x,y)\in S$ such that $y\le x/2$.
Then since $S_1$, $S_2$ and $S_3$ tile the set $\{(x,y)\in S:y\le x\}$,
we have
$$
E(\QQ^2)=2\int_{S_1\cup S_2\cup S_3} F\,dx\,dy\,.
$$
To estimate $F$ on $S_1$, we use the bound
\begin{align*}
P(x,y\in X)
&\le
P(x\in X)\wedge
P(y\in X)
\overset{\eqref{e.m1}}=
\frac{\yy(x)^\sb}{x^\sb}\wedge
\frac{\yy(y)^\sb}{y^\sb}
\,.
\end{align*}
Since $x/2\le y\le x$ on $S_1$, this is bounded by $\yy(y)^\sb\,y^{-\sb}\lesssim
\yy(y)^\sb\,x^{-\sb}$.
Hence,
\begin{multline*}
\int_{S_1} F(x,y)
\lesssim
\int_{\aa}^{\bb}
\int_{x-\yy(x)}^{x}
\frac{\rho(x)\rho(y)}{\yy(x)^\sb x^\sb}
\, dy\,dx
\overset{\eqref{e.regularity}}
\lesssim
\int_{\aa}^{\bb}
\int_{x-\yy(x)}^{x}
\frac{\rho(x)^2}{\yy(x)^\sb x^\sb}
\, dy\,dx
\\
=
\int_{\aa}^{\bb}
\frac{\rho(x)^2\,\yy(x)^{1-\sb}}{x^\sb}
\, dx
\end{multline*}
By the definition of $\Hk$, we have $\rho(x)\,\yy(x)^{1-\sb}\lesssim 1$.
Thus,~\eqref{e.Rp} implies that $\int_{S_1} F\lesssim 1$.

For $S_2$, we use the estimate~\eqref{e.m2}, the fact that $y\eqsim x$ when
$y\in [x/2,x-\yy(x)]$ and~\eqref{e.regularity}, to get
\begin{align*}
\int_{S_2} F
&
\lesssim \int_{\aa}^{\bb}\int_{x/2}^{x-\yy(x)}
\frac{\rho(x)}{x^\sb}\,\frac{\rho(y)}{(x-y)^\sb}\,dy\,dx
\\ &
\lesssim
\int_{\aa}^{\bb}\frac{\rho(x)^2}{x^\sb}\int_{x/2}^{x-\yy(x)}
\frac{1}{(x-y)^\sb}\,dy\,dx
\\ &
\lesssim
\int_{\aa}^{\bb}\frac{\rho(x)^2}{x^\sb}\,\Hk^{h}(x)^{-1}\,dx
=
\int_{\aa}^{\bb}\frac{\rho(x)}{x^\sb}\,dx=1\,.
\end{align*}
On the set $S_3$, the estimate~\eqref{e.m2} gives
$$
P(x,y\in X)\lesssim
\frac{\yy(x)^\sb \yy(y)^\sb}{x^\sb y^\sb}\,.
$$
Hence,
\begin{align*}
\int_{S_3} F
&
\lesssim
\int_{\aa}^{\bb}
\int_{\aa}^{\bb}
\frac{\rho(x)\rho(y)}{x^\sb y^\sb}\,dx\,dy
=
\Bigl(
\int_{\aa}^{\bb}
\frac{\rho(x)}{x^\sb}\,dx\Bigr)^2\overset{\eqref{e.m1}}=1\,.
\end{align*}
Thus, we conclude that
$ E(\QQ^2)\lesssim 1 = E(\QQ)^2$.
The Paley-Zygmund inequality therefore gives
\begin{equation}
\label{e.QQ}
P\bigl(\QQ\ge E\QQ/2\bigr)\gtrsim 1\,.
\end{equation}
Note that $\QQ>0$ implies $\sup X\ge \aa$.
But since $\aa$ can be arbitrarily large and the
constant implied in~\eqref{e.QQ} does not depend on $\aa$,
it follows from~\eqref{e.QQ} that $P(\sup X=\infty)\gtrsim 1$.

\medskip

Now fix some $t\in(0,\infty)$.
We will show that a.s.\
$$
P(\sup X=\infty \mid \mathcal F_t)\gtrsim 1\,.
$$
Define $X_t:=\bigl\{x\ge r: \tau^{h(x)}_x\le t\bigr\}$.
Suppose that $\sup X_t<\infty$.
Let $x>y>r\vee \sup X_t$.
Then we have by the
Markov property of $\SLE$ and~\eqref{e.m1} that
\begin{equation}
\label{e.Xt}
P(x\in X\mid \mathcal F_t)
=\bigl(\frac{g_t'(x)\,\yy(x)}{g_t(x)-W_t}\bigr)^{\sb }\,.
\end{equation}
The same reasoning, but this time with~\eqref{e.m2},
shows that
\begin{equation}
\label{e.Xt2}
P(x,y\in X \mid \mathcal F_t)
\lesssim
\bigl(\frac{g_t'(x)g_t'(y)\yy(x)\yy(y)}{(g_t(y)-W_t)(g_t(x)-g_t(y)}
\bigr)^\sb.
\end{equation}
Since $g_t$ has a power series expansion near $\infty$ of the form
$$
g_t(z)=z+\frac{a_1}z+\frac{a_2}{z^2}+\dots\,,
$$
we have in particular that $\lim_{x\to\infty} g_t(x)-x=0$
and $\lim_{x\to\infty} g'_t(x)=1$.
Therefore, on
the event $\sup X_t<\infty$
 there is some random $\aa'>\aa\vee \sup X_t$,
which is $\mathcal F_t$-measurable,
such that for all $x\ge \aa'$ we have $g_t'(x)\in[1/2,2]$
and
$(g_t(x)-W_t)/x \in[1/2,2]$.
Therefore, for $x>y>\aa'$ we have that the estimates in~\eqref{e.Xt}
and~\eqref{e.Xt2} are within a constant multiplicative factor
(which depends only on $\kappa$) from their values at $t=0$.
Consequently, our proof above with $\aa$ replaced by $\aa'$ and
with probabilities and expectations replaced by conditional
probabilities and conditional expectations given $\mathcal F_t$
implies that
on the event $\sup X_t<\infty$,
\begin{equation}
\label{e.supX}
P(\sup X>\aa'\mid \mathcal F_t)\gtrsim 1\,.
\end{equation}
But since $X\supset X_t$, this holds even if $\sup X_t=\infty$.
Because~\eqref{e.supX} a.s.\ holds for every $t>0$, we get $P(\sup X>\aa)=1$,
and since $\aa$ was arbitrary, we get $\sup X=\infty$ a.s.
Lemma~\ref{l.dist} implies therefore that the set of $x\ge r$
such that $\inf_t \dd^x_t \le 4\,\yy(x)$ is a.s.\ unbounded.

\medskip

Condition (\ref{e.regularity}) implies that for some finite constant $A>1$ and $x,y$ satisfying
$r\le x\le y\le 2\,x$, we have
$$
A\,h(y)\ge h(x) \qquad\text{if } \kappa<4\,,
$$
and
$$
\frac{h(y)}y \ge
\bigl(
\frac{ h(x)}x\bigr)^A \qquad\text{if } \kappa=4\,.
$$
%
Define
$$
H(x):=\begin{cases}(4A)^{-1} h(2x/3) 
&\kappa<4\,,\\
             4^{-1}(2x/3)^{1-A}h(2x/3)^A&\kappa=4\,.
\end{cases}
$$
Since the function $H(x)$ satisfies the same assumptions
we have for $\yy(x)$, we conclude that a.s.
$$
\sup\bigl\{x\ge 3\,r/2: \inf_t\dd^x_t\le 4H(x)\bigr\}=\infty\,.
$$
As the balls $B\bigl(x,4H(x)\bigr)$ lie below the graph of $y=h(x)$ when $x>2\,r$,
it follows that $\gamma\cap\{x+i\,y:y\le h(x),\,x\ge r\}$
is a.s.\ unbounded.
As we have seen, this implies that $\gamma\cap\Gamma^h$ is
a.s.\ unbounded.
The proof is thus complete.
\QED

\section{Hausdorff dimension} \label{HD}
In this part, we will prove Theorem~\ref{theore3}.
The usual strategy of deriving Theorem~\ref{theore3} is to estimate the two probabilities
$$
P(x\in \gamma_{\epsilon}  ), \
P(x\in \gamma_{\epsilon}, \ y\in \gamma_{\epsilon}),
$$
where
$\gamma_{\epsilon}:=\{x\in \R: \text{dist}(x,\gamma) \leq \epsilon\}$,
and then prove some $0$--$1$ law to show that the Hausdorff dimension is
an a.s.\ constant.

In this paper, instead of $\gamma_{\epsilon}$, we consider $C_{\epsilon}$.
Let
$$
C:=\bigcap_{\eps>0} C_\eps\,.
$$
Then Lemma~\ref{l.dist} gives
\begin{equation}
\label{e.Cg}
C\subset\gamma\cap\R\,.
\end{equation}

\begin{prop}\label{t2}
Assume that $\kappa\in(4,8)$.  Then for any $\delta>0$,
$$P(\dim_H C\geq 1-\sb-\delta )>0\,.$$
\end{prop}
\proof
The proof follows the standard Frostman measure argument.
We introduce random measures $\mu_{\epsilon}$
defined on the Borel $\sigma$-field of the interval $[1,2]$ by
$$
\mu_{\epsilon}([1,x]):=\epsilon^{-\sb}\int_1^x I(x_1\in C_{\epsilon})\,dx_1
$$
for $0<\epsilon<1$ and $x\in[1,2]$. The $(1-\sb-\delta)$-energy of $\mu_{\epsilon}$ is
$$
\mathcal{E}(\mu_{\epsilon})=\int_1^2\int_1^2 \f{1}{|y-x|^{1-\sb-\delta}}\,
d\mu_{\epsilon}(x)\,d\mu_{\epsilon}(y)\,.
$$
Its expectation is
\begin{eqnarray}
E\big(\mathcal{E}(\mu_{\epsilon})\big)&=&2\,\epsilon^{-2\sb}\int_1^2\int_x^2
\f{P(x\in C_{\epsilon},y\in C_{\epsilon}) }{|y-x|^{1-\sb-\delta}}\,dy\,dx \nonumber \\
&\leq& 2\,\epsilon^{-2\sb}\int_1^2\int_x^{x+\epsilon}
\f{P(x\in C_{\epsilon}) }{|y-x|^{1-\sb-\delta}}\,dy\,dx +{} \nonumber \\
& & \qquad {}+2\,\epsilon^{-2\sb}\int_1^{2-\epsilon}\int_{x+\epsilon}^2
\f{P(x\in C_{\epsilon},y\in C_{\epsilon}) }{|y-x|^{1-\sb-\delta}}\,dy\,dx \nonumber \\
&=:&E_1+E_2\,. \label{E12}
\end{eqnarray}
For $E_1$,  Proposition~\ref{m1} gives that
\begin{equation}
E_1 \leq 2\,\epsilon^{-\sb}\int_1^2\int_x^{x+\epsilon}(y-x)^{\sb-1}\,dy\,dx 
=2\,\sb^{-1}. \label{E1}
\end{equation}
For $E_2$, Proposition~\ref{m2} gives that
\begin{equation}
E_2\leq 2\,c_{\kappa}\int_1^{2-\eps}\int_{x+\epsilon}^2(y-x)^{-1+\delta}\,dy\,dx 
\leq  2\,c_{\kappa}\,\delta^{-1}. \label{E2}
\end{equation}
Combining (\ref{E12}), (\ref{E1}) and (\ref{E2}), we obtain that
$$
E\big(\mathcal{E}(\mu_{\epsilon})\big) \leq 2\,\sb^{-1}+2\,c_{\kappa}\,\delta^{-1}.
$$
Noting that
\begin{equation} \label{excon}
E|\mu_{\epsilon}|=\epsilon^{-\sb}\int_1^2 (\epsilon/x)^{\sb}\,
dx=(1-\sb)^{-1}(2^{1-\sb}-1)>0\,,
\end{equation}
and $E\bigl(|\mu_\eps|^2\bigr)\le E\bigl( \mathcal E(\mu_\eps)\bigr)$,
the Paley-Zygmund inequality implies
that there is a $\lambda>0$, which does not depend on $\eps$, such that with probability
at least $\lambda$,
$|\mu_{\epsilon}|>\lambda$ and $\mathcal{E}(\mu_{\epsilon})<1/\lambda$.
With probability at least $\lambda$ this will hold for a sequence of positive $\eps$
tending to $0$.
On this event, we can take a subsequential limit $\mu$
supported on $C$ and satisfying
$|\mu|>\lambda$ and $\mathcal{E}(\mu)<1/\lambda$.
Frostman's lemma therefore implies that $P\bigl(\dim_H(C\cap[1,2])>1-\sb-\delta\bigr)>\lambda$,
which concludes the proof.
\QED

The following proposition tells us that $\dim_H\bigl(\gamma\cap\R\bigr)\le 1-\sb $ a.s.

\begin{prop}\label{p1}
Let $x\in [1,2)$, and $\kappa\in(4,8)$. Then for $\eps\in(0,1)$,
$$P(\gamma \cap [x,x+\epsilon] \neq \emptyset)
\asymp  \left(\frac{\epsilon}{x}\right)^{\sb },$$
where the constants implied by $\asymp$ depend only on $\kappa$.
\end{prop}
\proof Lemma 6.6 of \cite{RS} gives
\begin{multline}
\label{6.6}
P(\gamma \cap [x,x+\epsilon] \neq \emptyset) \\
=1-\frac{4^{(\kappa-4)/\kappa}\sqrt{\pi}\,
 {}_2F_1(1-4/\kappa,2-8/\kappa,2-4/\kappa;1/q)
q^{(4-\kappa)/\kappa}}{\Gamma(2-4/\kappa)\Gamma(4/\kappa-1/2)},
\end{multline}
where $q:=(x+\epsilon)/x$. Using
$\Gamma(2\theta)\,\Gamma(1/2)=2^{2\theta-1}\,\Gamma(\theta)\,\Gamma(\theta+1/2)$,
 Euler's integral representation of ${}_2F_1$ and a change of variable, we have
\begin{eqnarray*}
(\ref{6.6})&=&1-\frac{\Gamma(4/\kappa)}
{\Gamma(1-4/\kappa)\Gamma(8/\kappa-1)} \, \int_0^{\frac{x}{x+\epsilon}}
t^{-4/\kappa} (1-t)^{8/\kappa-2} \; dt\notag\\
&=&\frac{\Gamma(4/\kappa)}
{\Gamma(1-4/\kappa)\Gamma(8/\kappa-1)} \, \int_{\frac{x}{x+\epsilon}}^1
t^{-4/\kappa} (1-t)^{8/\kappa-2} \; dt\notag\\
&\asymp& \left(\frac{\epsilon}{x}\right)^{\sb }.
\QED
\end{eqnarray*}

\begin{lemma} \label{zeroone}
There is a constant $d=d_\kappa$ such that
$\dim_H (\gamma \cap \R)=d$ a.s.
\end{lemma}
\proof For all $n\in \Z$ let $D_n:=\dim_H\bigl( \gamma[0,2^n]\cap \R\bigr)$.
Then $D_{n+1}\geq D_n$.
In addition, $D_n$ and $D_{n+1}$ have the same distribution, by scale invariance.
Therefore, $D_m=D_n$ a.s.\ for all $m,n\in\Z$. Hence,
$\dim_H( \gamma\cap \R)=\sup_{n\in\Z}D_n$ is
$\mathcal{F}_{2^n}$-measurable for all $n\in\Z$, which implies that $\dim_H( \gamma\cap \R)$
 is $\mathcal{F}_{0^+}$-measurable.
 By Blumenthal's $0$-$1$ law, the $\sigma$-field $\mathcal F_{0^+}$ is trivial. \QED

\proofof{Theorem~\ref{theore3}}
First note that $\dim_H (\gamma\cap \R)=\dim_H(\gamma\cap \R_+)=\dim_H(\gamma\cap \R_-)$ a.s.\ by the symmetry
property of $\SLE$ curves.
 Proposition~\ref{p1} implies $\dim_H\bigl( \gamma\cap \R^+\bigr)
\leq 1-\sb$  a.s. On the other hand,~\eqref{e.Cg}  and
Proposition~\ref{t2} give
$$P\bigl(\dim_H (\gamma \cap \R) \geq 1-\sb-\delta \bigr)>0$$
for every $\delta>0$.
Therefore $\dim_H( \gamma \cap \R)=1-\sb $ a.s.\ by
Lemma~\ref{zeroone}. \QED

\bigskip\noindent{\bf Acknowledgment}.
We are grateful to David Wilson and Yuval Peres for helpful discussions.
We thank Scott Sheffield for comments and suggestions regarding
 a previous version of this paper.


\begin{thebibliography}{00}

\bibitem{AS}
T.~Alberts and S.~Sheffield. (2007).
Hausdorff Dimension of the $\SLE$ curve intersected with
the real line. \mbox{arXiv:0711.4070v1}.

\bibitem{beffara04}
V.~Beffara. (2004).
Hausdorff dimensions for $\SLE$(6).
{\em Ann. Prob.}, 32:2606--2629.

\bibitem{beffara}
V.~Beffara. (2007).
 The dimension of the $\SLE$ curves.
{\em To appear in Ann. Prob.}.

\bibitem{vdB-Jarai}
R.~van den Berg, A.~A.~J\'arai. (2003).
The lowest crossing in two-dimensional critical percolation.
{\em Ann. Prob.}, 31:1241--1253.

\bibitem{CN}
F.~Camia and C.~M. Newman. (2005).
The Full Scaling Limit of Two-Dimensional
  Critical Percolation. \mbox{arXiv:math.Pr/0504036}.

\bibitem{cardy}
J.~Cardy. (2005).
$\SLE$ for theoretical physicists.
{\em Ann. Physics}, 318:81--118.

\bibitem{GK}
 I.~A. Gruzberg and L.~P. Kadanoff. (2004).
 The Loewner equation: maps and shapes.
 {\em J. Statist. Phys.}, 114:1183--1198.

\bibitem{lawler}
G.~F.~Lawler. (2005).
 {\em Conformally Invariant Processes in the Plane}.
 American Mathematical Society, Providence, RI.

\bibitem{LSW04}
G.~Lawler, O.~Schramm and W.~Werner. (2004).
Conformal invariance of planar loop-erased random walks and uniform
spanning trees.
{\em Ann. Prob.}, 32:939--995.


\bibitem{RS}
S.~Rohde and O.~Schramm. (2005).
 Basic properties of $\SLE$.
 {\em Ann. Math.}, 161:883--924.

\bibitem{sch00}
 O.~Schramm. (2000).
  Scaling limits of loop-erased random walks and uniform
  spanning trees. {\em Israel Journal of Mathematics}, 118:221--288.


\bibitem{SS08}
O.~Schramm and S.~Sheffield. (2008).
Contour lines of the two-dimensional discrete Gaussian free field.
  {\em Acta Math.}, to appear.
 \mbox{arXiv:math/0605337}.

\bibitem{SW}
O.~Schramm and D. B.~Wilson. (2005).
$\SLE$ coordinate changes.
 {\em New York Journal of Mathematics}, 11:659--669.

\bibitem{smirnov01}
 S.~Smirnov. (2001).
 Critical percolation in the plane: Conformal invariance,
  Cardy's formula, scaling limits.
  {\em C. R. Acad. Sci. Paris S\'er. I Math.},
  333:239--244.

\bibitem{smirnov07}
 S.~Smirnov. (2007).
 Conformal invariance in random cluster models.
 I. Holomorphic fermions in the Ising model.
 \mbox{arXiv:0708.0039v1}.

\bibitem{werner}
W.~Werner. (2004). \emph{Random planar curves and Schramm-Loewner evolutions},
  Lectures on probability theory and statistics,
  Lecture Notes in Math.,
  1840, pp.~107--195,  Springer, Berlin.  \mbox{arXiv:math.PR/0303354}.

\bibitem{zhan07}
D.~Zhan. (2007).
 Reversibility of chordal SLE.
 \mbox{arXiv:0705.1852}.
\end{thebibliography}
\end{document}